\documentclass[3p,sort&compress]{elsarticle}
\usepackage{amsmath,amssymb,amsthm}
\usepackage{graphicx}

\begin{document}

\journal{{}}

\begin{frontmatter}
\title{Numeric Solution of Advection-Diffusion Equations by a Discrete Time Random Walk Scheme} 

\author[UNSW]{C. N. Angstmann\corref{cor1}}
\ead{c.angstmann@unsw.edu.au}
\author[UNSW]{B. I. Henry}
\author[WITS,COE]{B. A. Jacobs}
\author[UNSW]{A. V. McGann}

\address[UNSW]{School of Mathematics and Statistics, UNSW Australia, Sydney NSW 2052 Australia}
\address[WITS]{
	School of Computer Science and Applied Mathematics, 
			University of the Witswatersrand, 
			Johannesburg, 
			Private Bag 3, 
			Wits 2050, 
			South Africa}
\address[COE]{
			DST-NRF Centre of Excellence in Mathematical and Statistical Sciences (CoE-MaSS),
            University of the Witwatersrand, 
            Johannesburg,
            Private Bag 3, 
            Wits 2050, 
            South Africa}
\cortext[cor1]{Corresponding author}

\begin{abstract}
Explicit numerical finite difference schemes for partial differential equations are well known to be easy to implement
but they are particularly problematic for solving equations
whose solutions admit shocks, blowups and discontinuities.  
Here we present an explicit numerical scheme for solving non-linear advection-diffusion equations
admitting shock solutions
that is both easy to implement and stable. 
The numerical scheme is obtained by considering the continuum limit of a discrete time and space stochastic process for non-linear advection-diffusion. The  stochastic process is well posed and this guarantees the stability of the scheme.  Several examples are provided to highlight the importance of the
formulation of the stochastic process in obtaining a stable and accurate numerical scheme. \end{abstract}

\end{frontmatter}

\section{Introduction}

This work considers a general form of nonlinear dissipative dominated partial differential equations and more closely examines the numerical solutions of the Burgers' equation which is a limiting form of the more general dissipative system, as derived by Su and Gardner \cite{su1969korteweg}. The Burgers' equation admits some difficulties in obtaining numerical solutions under sufficiently high Reynolds numbers wherein shocks in the solution may form. Many authors have obtained numerical solutions to the Burgers' equation using finite elements, finite differences, exponentially fitted methods, spectral methods, cubic B-splines, Adomian decomposition methods, differential quadrature, wavelets, compact schemes and method of lines  \cite{caldwell1982solution,kutluay1999numerical,hon1998efficient,ozics2009exponentially,mlttal1996numerical,daug2005numerical,abbasbandy2005numerical,bellman1972differential,beylkin1997adaptive,bhatt2016fourth,mukundan2015efficient}.

Recently Klein and Saut \cite{klein2015numerical} introduce a novel approach to the Burgers' equation as well as the fractional Korteweg-de Vries equations. Angstmann et al. \cite{ADH2013mmnp,ADHN2014,angstmann2015generalized,angstmann2016stochastic} provide a framework for the construction of fractional and integer order reaction-diffusion models which is able to recover a general class of parabolic-hyperbolic equations. This framework is based on the random walk principles. In \cite{angstmann2016stochastic} Angstmann et al. illustrate how the discrete time random walk is able to numerically solve equations of this class.

Thom{\'e}e \cite{thomee2001finite} provides an extensive review of finite difference and finite element analysis for partial differential equations. This work covers the concepts of convergence and stability for numerical methods for partial differential equations as well as the importance of the Courant, Freidrichs and Lewy (CFL) condition for numerical solutions to advective equations.

Exact solutions to the Burgers' equation have been sought through the well known Hopf-Cole \cite{Hopf1950,Cole1951} transformation, which recovers the standard heat equation from the Burgers' equations. Fletcher \cite{fletcher1983generating} generates exact solutions for the two-dimensional Burgers' equation through the multi-dimensional interpretation of the Hopf-Cole transformation.

This work derives a master equation from a discrete time and space stochastic process as a numerical scheme to approximate a partial differential equation. The correspondence between the discrete master equation and the partial differential equation is established by taking the diffusion limit of the discrete stochastic process. 

We first construct the master equation that will govern the evolution of the probability mass function of the random walking particle. We then show that, under the appropriate limit, this master equation will approach an advection diffusion equation. This shows that the master equation can function as an approximation of the advection diffusion equation. Further more, as the limit process of the random walk exists under the same limit, the solution of the difference equation will converge towards the solution of the advection diffusion equation. This allows us to construct a numerical scheme for the advection diffusion equation. 

Appropriate construction of the jump probabilities of the random walk guarantees that the stochastic process remains well defined for any given lattice spacing. Similar constructions are also used on the boundary conditions to again guarantee that the random walk remains valid for all lattice spacings.  

A rigorous numerical analysis of the present method is conducted as a means of verifying the results and illustrating the efficacy of the method.

\section{A Discrete Time Random Walk}

\subsection{Master Equation}
We consider a discrete time random walk on a one dimensional lattice. A particle begins at an initial position and at each time step the particle will jump from its current lattice site to a neighbouring lattice site. The probability that a particle on the $i^{\mathrm{th}}$ lattice site will jump to the right on the $n^{\mathrm{th}}$ time step is given by an arbitrary time and space dependent probability, $P_r(i,n)$, and the probability of jumping to the left is given by the complement, $P_{l}(i,n)=1-P_{r}(i,n)$. We can recursively define the probability of the particle being at lattice site $i$ at time step $n$, denoted $U(i,n)$, given some initial condition, $U(i,0)$, 
\begin{equation}
\label{eq_dme}
U(i,n)= P_{r}(i-1,n-1)U(i-1,n-1)+P_{l}(i+1,n-1)U(i+1,n-1).
\end{equation}
This is the master equation governing the time evolution of the probability mass function for the process. Here we have considered an infinite lattice, later we will discuss a finite lattice, in which boundary conditions $U(L,n)$ are defined for a boundary set at $x=L$ at all times $n$. 

\subsection{Diffusion Limit}

We wish to consider the continuum space and time limit of the discrete time random walk. The appropriate limit in this case is a diffusion limit in which the space and time scales are together taken to zero with the requirement that the following limit exists,
\begin{equation}
\label{eq:limdif}
D=\lim_{\Delta x, \Delta t \to 0}\frac{\Delta x^2}{2 \Delta t},
\end{equation}
where $\Delta x$ is the spatial grid spacing and $\Delta t$ is the temporal grid spacing \cite{Kac1947}. In the diffusion limit the master equation, Eq. (\ref{eq_dme}), will become an advection-diffusion PDE.  As such, the master equation itself may be considered as an approximation to the advection-diffusion equation.

First we begin by associating each of the discrete functions, $U(i,n)$ with a continuous function,$u_{\Delta}(x,t)$, such that,
\begin{equation}
u_{\Delta}(i\Delta x , n\Delta t)=U(i,n).
\end{equation}
This association is dependent on the lattice spacings, $\Delta x$ and $\Delta t$. The jump probabilities may also be associated with continuum functions, sampled at discrete points,  
\begin{equation}
p_{r,\Delta}(i\Delta x,n \Delta t)=P_r(i,n),
\end{equation}
and similarly for $P_l(i,n)$. At each point the continuum jump probabilities sum to one,
\begin{equation}
p_{r,\Delta}(i\Delta x,n \Delta t)+p_{l,\Delta}(i\Delta x,n \Delta t)=1.
\end{equation}
 At the point $x=i \Delta x$ and $t=n \Delta t$ ,the discrete master equation, Eq. (\ref{eq_dme}) then becomes,
 \begin{equation}
u_{\Delta}(x, t)= p_{r,\Delta}(x-\Delta x ,t- \Delta t)u_{\Delta}(x-\Delta x, t-\Delta t )+p_{l,\Delta}(x+\Delta x ,t- \Delta t)u_{\Delta}(x+\Delta x ,t- \Delta t ).
\end{equation}
Defining a force function, 
\begin{align}
f_{\Delta}(x,t)&=p_{r,\Delta}(x,t)-p_{l,\Delta}(x,t),\\
&=2p_{r,\Delta}(x,t)-1,
\end{align}
the master equation can be expressed as,
\begin{equation}
\begin{split}
u_{\Delta}(x, t)=&\frac{1}{2}\left(u_{\Delta}(x-\Delta x, t-\Delta t )+u_{\Delta}(x+\Delta x ,t- \Delta t )\right)\\
&+\frac{f_{\Delta}(x-\Delta x,t-\Delta t)}{2}\left(u_{\Delta}(x-\Delta x, t-\Delta t )\right)-\frac{f_{\Delta}(x+\Delta x,t-\Delta t)}{2}\left(u_{\Delta}(x+\Delta x, t-\Delta t )\right)
\end{split}
\end{equation}
Further manipulations give,
\begin{equation}
\begin{split}
\label{eq_nonlim}
&\frac{u_{\Delta}(x, t)-u_{\Delta}(x, t-\Delta t)}{\Delta t}=\frac{\Delta x^2}{2 \Delta t}\left(\frac{u_{\Delta}(x-\Delta x, t-\Delta t )-2u_{\Delta}(x, t-\Delta t)+u_{\Delta}(x+\Delta x ,t- \Delta t )}{\Delta x^2}\right)\\
&+\frac{\Delta x^2}{2\Delta t}\left(\frac{f_{\Delta}(x-\Delta x,t-\Delta t)\left(u_{\Delta}(x-\Delta x, t-\Delta t )\right)-f_{\Delta}(x+\Delta x,t-\Delta t)\left(u_{\Delta}(x+\Delta x, t-\Delta t )\right)}{\Delta x^2}\right)
\end{split}
\end{equation}

We are now set up to take the diffusion limit, that is the limit $\Delta x \to 0$ and $\Delta t \to 0$ such that the limit, Eq. (\ref{eq:limdif}) exists.
This limit will take the random walk stochastic process to a limit process provided that the limit process exists. In this case the limit process will exist as the original stochastic process is a simple random walk \cite{GS2016}. The particle distribution from the limit process will be given by the limit of the continuum embedding of the discrete process distribution, i.e.,
\begin{equation}
u(x,t)=\lim_{\Delta x, \Delta t \to 0} u_{\Delta}(x,t).
\end{equation}
Taking the limit of Eq. (\ref{eq_nonlim}) gives,
\begin{equation}
\begin{split}
\label{eq_cl}
\frac{\partial u(x, t)}{\partial t}=&D \frac{\partial ^{2}u(x, t)}{\partial x^2}- 2\beta D \frac{\partial}{\partial x}\left(F(x,t)u(x,t)\right),
\end{split}
\end{equation}
where,
\begin{equation}
\beta F(x,t)=\lim_{\Delta x\to 0}\frac{f_{\Delta}(x,t)}{\Delta x}=\lim_{\Delta x\to 0}\frac{p_{r,\Delta}(x,t)-p_{l,\Delta}(x,t)}{\Delta x}.
\end{equation}

This shows that the diffusion limit of the random walk master equation is an advection-diffusion PDE. It should be noted that there is no condition on the form of the force function and this PDE may be a non-linear advection-diffusion equation by considering $F(x,t)$ as an explicit function of $u(x,t)$.

\section{Formulation of the Numerical Scheme}

Given a general non-linear advection diffusion PDE of the form 
\begin{equation}
\label{eq_genAD}
\frac{\partial u(x,t)}{\partial t}=D\frac{\partial^2 u(x,t)}{\partial x^2}-2\beta D \frac{\partial}{\partial x}\left(F(x,t)u(x,t)\right),
\end{equation}
we can form a numerical scheme by matching the diffusion limit of the master equation to the PDE.
The diffusion limit of the master equation is completely defined by two parameters, $\Delta x$, and $\Delta t$, as well as a single function $P_r(i,n)$, and hence we need a method of obtaining these parameters from the PDE. Firstly, fixing the spatial lattice spacing $\Delta x$ will also fix the temporal step size through the relation,
\begin{equation}
\Delta t=\frac{\Delta x^2}{2 D}.
\end{equation}
Next, we have a condition on the function $P_r(i,n)$ that is imposed by matching the PDE to the diffusion limit of the master equation. Given a $\Delta x$ and $\Delta t$ we have $p_{r,\Delta}(i \Delta x,n \Delta t)=P_r(i,n)$ and in the diffusion limit we require,
\begin{equation}
\label{eq_f_lim}
F(x,t)=\frac{f_\Delta(x,t)}{\beta}=\lim_{\Delta x, \Delta t \to 0}\frac{2 p_{r,\Delta}(x,t)-1}{\beta \Delta x}.
\end{equation}
Whilst any choice of function, $p_{r,\Delta} (x,t)$, that obeys this limit will work it is desirable to chose a function that will provide a stable numerical scheme far from the limit. In order for the master equation to describe a valid stochastic process the jump probability must always lie between $0$ and $1$, i.e. $0\leq P_r(i,n)\leq 1$. The most obvious form from the limit relation, Eq. (\ref{eq_f_lim}), would be to take, 
\begin{equation}
 p_{r,\Delta}(x,t)=\frac{\beta F(x,t) \Delta x+1}{2}.
\end{equation}
This is a poor choice as the restriction $p_{r,\Delta}(x,t)\leq1$, puts an upper bound on $\Delta x$.  To overcome this limitation we construct the jump probability from Boltzmann weights \cite{HLS2010}. These weights are taken by considering a diffusing particle at equilibrium in a potential $V(x,t)$ where,
\begin{equation}
\label{eq_f_pot}
F(x,t)=-\frac{\partial V(x,t)}{\partial x}.
\end{equation}
The position of such a particle will follow a Boltzmann distribution such that the probability of finding the particle at position $x$ at time $t$ will be proportional to $\exp(-\beta V(x,t))$, where $\beta$ is related to the ``temperature" of the system. If we restrict this distribution such that the particle may only be present at $x+\Delta x$ or $ x-\Delta x$ then the normalisation constant is simply  $\exp(-\beta V(x-\Delta x,t))+\exp(-\beta V(x+\Delta x,t))$. Given that the particle begins at $x$, we interpret the probability that the particle has jumped right from $x$ is the Boltzmann probability of the particle being at $x+\Delta x$. Similarly we interpret the probability that the particle has jumped left from $x$ as the Boltzmann probability of the particle being at $x-\Delta x$.We can then write,
\begin{equation}
\label{eq_boltz}
 p_{r,\Delta}(x,t)=\frac{\exp(-\beta V(x+\Delta x,t))}{\exp(-\beta V(x-\Delta x,t))+\exp(-\beta V(x+\Delta x,t))}.
\end{equation}
Such a functional form will obey the limit condition, Eq. (\ref{eq_f_lim}), as well as guarantee that $0\leq  p_{r,\Delta}(x,t)\leq1$ for all $\Delta x$ and $F(x,t)$. Using this jump probability the master equation can be written with an arbitrary $\Delta x$ and remain a valid stochastic process. 

The Boltzmann weights are defined in terms of a potential $V(x,t)$. In order to implement the numerical scheme we need to rework this into a function of the force $F(x,t)$. From Eq. (\ref{eq_f_pot}) we can write,
\begin{equation}
V(x,t)=-\int_{0}^{x}F(x',t)dx'.
\end{equation}
The Boltzmann weights call for the potential evaluated at $x+\Delta x$ and $x-\Delta x$, these can be written relative to the potential at $x$,
\begin{equation}
\begin{split}
V(x+\Delta x,t)&=-\int_{0}^{x}F(x',t)dx'-\int_{x}^{x+\Delta x}F(x',t)dx'\\
&=V(x,t)-\int_{x}^{x+\Delta x}F(x',t)dx',
\end{split}
\end{equation}
and similarly,
\begin{equation}
\begin{split}
V(x-\Delta x,t)&=-\int_{0}^{x}F(x',t)dx'+\int_{x-\Delta x}^{x}F(x',t)dx'\\
&=V(x,t)+\int_{x-\Delta x}^{x}F(x',t)dx'.
\end{split}
\end{equation}
If $F(x,t)$ is a simple function of space and time then this integral may be solved analytically and the results used in the jump probabilities. A complication arrises when the force is dependent on the particle distribution itself, i.e. a non linear force $F(x,t,u(x,t))$. While the above integral is conceptually fine, we will only have an approximation to the probability distribution on the lattice points. In computing the integral we will need to use a numerical quadrature that is restricted to these lattice points.
There are two convenient approximations that can be used, either a single point or a two point quadrature.  The single point quadrature is,
\begin{align}
\int_{x}^{x+\Delta x}F(x',t)dx'\approx \Delta x F(x,t)\\
\int_{x-\Delta x}^{x}F(x',t)dx'\approx \Delta x F(x,t)
\end{align}
Substituting the resulting approximations for the potential into the Boltzmann weights, Eq. (\ref{eq_boltz}), and simplifying gives,
\begin{equation}
\label{eq_boltz1p}
 p_{r,\Delta}(x,t)=\frac{\exp(\beta\Delta x F(x,t))}{\exp(-\beta\Delta x F(x,t))+\exp(\beta\Delta xF(x,t))}.
\end{equation}
or more compactly,
\begin{equation}
\label{eq_boltz1pa}
 p_{r,\Delta}(x,t)=\frac{1}{1+\exp(-2\beta\Delta x F(x,t))}.
\end{equation}
The single point quadrature results in jump probabilities that only require the force to be evaluated at a single lattice point. 

The two point quadrature is,
\begin{align}
\int_{x}^{x+\Delta x}F(x',t)dx'\approx\frac{\Delta x}{2}\left(F(x,t)+F(x+\Delta x,t)\right)\\
\int_{x-\Delta x}^{x}F(x',t)dx'\approx \frac{\Delta x}{2}\left(F(x-\Delta x,t)+F(x,t)\right)
\end{align}
Substituting the resulting approximations for the potential into the Boltzmann weights, Eq. (\ref{eq_boltz}), and simplifying gives,
\begin{equation}
\label{eq_boltz2}
 p_{r,\Delta}(x,t)=\frac{\exp(\frac{\beta\Delta x}{2}\left( F(x+\Delta x,t)+F(x,t)\right))}{\exp(-\frac{\beta\Delta x}{2} \left(F(x-\Delta x,t)+F(x,t)\right))+\exp(\frac{\beta\Delta x}{2} \left(F(x+\Delta x,t)+F(x,t)\right))}.
\end{equation}
or more compactly,
\begin{equation}
\label{eq_boltz3}
 p_{r,\Delta}(x,t)=\frac{1}{1+\exp(-\frac{\beta\Delta x}{2} \left(F(x-\Delta x,t)+2F(x,t)+F(x+\Delta x,t)\right))}.
\end{equation}
To evaluate the jump probabilities using the two point quadrature requires evaluating the force at three lattice points.

The final discrete time random walk numerical scheme for a non-linear advection diffusion equation can then be found by substituting these weights into the master equation. The two point quadrature jump probabilities give,
\begin{equation}
\label{eq_NS}
\begin{split}
u(i\Delta x,n\Delta t)= &\frac{u((i-1)\Delta x,(n-1)\Delta t)}{1+\exp(-\frac{\beta\Delta x}{2} \left(F((i-2)\Delta x,(n-1)\Delta t)+2F((i-1)\Delta x,(n-1)\Delta t)+F(i\Delta x,(n-1)\Delta t)\right))}\\
&+\frac{u((i+1)\Delta x,(n-1)\Delta t)}{1+\exp(\frac{\beta\Delta x}{2} \left(F((i+2)\Delta x,(n-1)\Delta t)+2F((i+1)\Delta x,(n-1)\Delta t)+F(i\Delta x,(n-1)\Delta t)\right))}.
\end{split}
\end{equation}
Using the two point quadrature jump probabilities results in a method that requires knowledge of the function at five lattice points. Using the single point quadrature results in a method that only requires three points. 

In this formulation we have a single free scale parameter, which may be cast as either $\Delta x$ or $\Delta t$. Every other parameter or function comes from the PDE that we wish to approximate. 

\subsection{Initial Conditions}
\label{sec_ic}
The scheme has been derived from the evolution of a probability density function. In general the integral of the initial condition over the domain will not equal one, and hence cannot be a probability density. As such we note that it is always possible to rescale any initial condition and we can consider the evolution equation to evolve a general non-normalised distribution. 
Another consequence of the distributional nature of the derivation is that we have a guarantee that the solutions will conserve mass, provided that the boundary conditions preserve mass, and be bound non-negative. Whilst this is often advantageous, a further complication arises if the initial condition is not strictly non-negative. This can not be interpreted as a distribution and hence may not be evolved forward by the master equation. However it is possible to write any generalised function as the difference between two distributions. We simply define two strictly non-negative distributions $u^{+}(x,t)$ and $u^{-}(x,t)$ such that,
\begin{equation}
u(x,t)=u^{+}(x,t)-u^{-}(x,t).
\end{equation}
If the strictly non-negative distributions evolve according to the coupled non-linear advection diffusion equations,
\begin{align}
\frac{\partial u^{+}(x,t)}{\partial t}&=\frac{\partial^2 u^{+}(x,t)}{\partial x^2}-\nu \frac{\partial F(x,t,u^{+}-u^{-}) u^{+}(x,t)}{\partial x},\\
\frac{\partial u^{-}(x,t)}{\partial t}&=\frac{\partial^2 u^{-}(x,t)}{\partial x^2}-\nu \frac{\partial F(x,t,u^{+}-u^{-}) u^{-}(x,t)}{\partial x},
\end{align}
then $u(x,t)$ will evolve according to the non-linear advection diffusion equation,
\begin{equation}
\frac{\partial u(x,t)}{\partial t}=\frac{\partial^2 u(x,t)}{\partial x^2}-\nu \frac{\partial F(x,t,u) u(x,t)}{\partial x}.
\end{equation}
Note that the coupling occurs through the non-linear force term. The evolution of the distributions $u^{+}$ and $u^{-}$ may be approximated as coupled DTRWs. The resulting numerical scheme is,
\begin{align}
\label{eq_NS2a}
U^{+}(i,n)= &\frac{U^{+}(i-1,n-1)}{1+\exp(-\frac{\beta\Delta x}{2}\left(F(i-2,n-1)+2F(i-1,n-1)+F(i,n-1)\right))}\\
&+\frac{U^{+}(i+1,n-1)}{1+\exp(\frac{\beta\Delta x}{2} \left(F(i+2,n-1)+2F(i+1,n-1)+F(i,n-1)\right))}. \nonumber\\
\label{eq_NS2b}
U^{-}(i,n)= &\frac{U^{-}(i-1,n-1)}{1+\exp(-\frac{\beta\Delta x}{2} \left(F(i-2,n-1)+2F(i-1,n-1)+F(i,n-1)\right))}\\
&+\frac{U^{-}(i+1,n-1)}{1+\exp(\frac{\beta\Delta x}{2}\left(F(i+2,n-1)+2F(i+1,n-1)+F(i,n-1)\right))}.\nonumber
\end{align}
where $U(i,n)=U^{+}(i,n)-U^{-}(i,n)$. Notice that if we take the difference of Eq. (\ref{eq_NS2a}) and Eq. (\ref{eq_NS2b}) then we recover the original numerical scheme, Eq. (\ref{eq_NS}) with no appearance of $U^{+}$ or $U^{-}$. Thus we can use the original numerical scheme for the case of initial conditions with mixed sign, provided that the boundary conditions also do not contain $U^{+}$ or $U^{-}$.

\subsection{Boundary Conditions}

Care must be taken with the implementation of boundary conditions so that the underlying random walk remains valid. The boundary conditions thus need to be consistent with both the underlying stochastic process and the given PDE boundary conditions. In general, this excludes some conditions, but does still allow for the implementation of a wide variety. The scheme can accommodate non-conservative boundary conditions, but care needs to be taken that the boundaries do not force a positive distribution to become negative. 

\subsubsection{Dirichlet Boundary Conditions}

In the case of a strictly non-negative initial condition, Dirichlet boundary conditions are also restricted to be non-negative. The conditions are implemented by fixing the value at the boundary. Using the two point quadrature, the jump probabilities are dependent on the solution at neighbouring points. As the jump probabilities need to be evaluated at the boundary point an additional point is required to use the two point quadrature form of the jump probabilities at the boundary. If an additional point is not known then the single point quadrature is required to be used.  

Given the boundary condition at $x=l$,
\begin{equation}
u(l,t)=a(t),
\end{equation}
where $a(t) \in \mathbb{R}^{+}$, we set 
\begin{equation}
\label{eq_DBC}
U(L,n)=a(n \Delta t)
\end{equation}
where $l=L\Delta x$. 

In the case that the initial condition is of mixed sign, then the Dirichlet boundary condition must be compatible with the extension given in Section \ref{sec_ic}. Again breaking the solution into two non-negative distributions, $u(x,t)=u^+(x,t)-u^-(x,t)$,  the boundary condition needs to be formulated for both non-negative distributions. Given a boundary condition at $x=l$ of the form,
\begin{equation}
u(l,t)=a(t),
\end{equation} 
we will require,
\begin{equation}
u^{+}(l,t)-u^{-}(l,t)=a(t).
\end{equation} 
Given that we also require $u^{+}(l,t)\geq0$ and $u^{-}(l,t)\geq0$, this is most easily satisfied by,
\begin{equation}
u^{+}(l,t)=\max(a(t),0),
\end{equation} 
and
\begin{equation}
u^{-}(l,t)=\max(-a(t),0).
\end{equation} 
The boundary conditions are then implemented by setting, 
\begin{equation}
U^{+}(L,n)=\max(a(n \Delta t),0),
\end{equation}
and 
\begin{equation}
U^{-}(L,n)=\max(-a(n \Delta t),0),
\end{equation}
where $l=L\Delta x$. 

\subsubsection{Neumann Boundary Conditions}

The restriction on Neumann boundary conditions are more subtle. We will consider here the case of a non-negative initial condition, although the results are again simply extended using the approach in Section \ref{sec_ic}. In the same manner as the standard finite difference schemes, the boundary condition may be implemented by including a ghost point outside the domain. Given a Neumann boundary condition at $x=l$ of the form,
\begin{equation}
\frac{\partial u(x,t)}{\partial x}\Big\vert_{x=l}=b(t)
\end{equation}
where $b(t) \in \mathbb{R}$, and where the domain sits to the left of $x=l$, we may set the ghost point via standard finite difference approximations,
\begin{equation}
\label{eq_NBCs}
U(L+1,n)=U(L,n)+\Delta x b(n \Delta t) 
\end{equation}
where $l=\left(L+\frac{1}{2}\right) \Delta x$. This will result in the correct behaviour of the numerical scheme as $\Delta x \to 0$, but may not be consistent with the stochastic process. If $\Delta x$ is sufficiently large, depending on the exact form of $b(t)$, it would be possible for the value at the ghost points to become negative. An alternative is to construct the ghost points by scaling the boundary points. In order for the boundary condition to hold we require,
\begin{equation}
\lim_{\Delta x, \Delta t \to 0} \frac{U(L+1,n)-U(L,n)}{\Delta x}=b(t), 
\end{equation}
where $t=n\Delta t$.
This limit relationship will be obeyed by the standard finite difference approximations, but there are many other choices that could be made. In keeping with the Boltzmann weights construction we may take,
\begin{equation}
\label{eq_NBC}
U(L+1,n)=U(L,n) \exp\left(\Delta x \frac{b(n \Delta t)}{U(L,n)}\right).
\end{equation}
These results are easily extendable for considering a left-boundary point, where the domain sits to it's right. This choice guarantees that the ghost points remain positive and hence that they remain consistent with the stochastic process. Once again the single point quadrature form of the jump probabilities needs to be used at the ghost point.

\section{Numerical Analysis}
	\subsection{Stability}
	
	If the distance between the approximation from the numerical scheme and the solution to the PDE is bounded for all time steps then the numerical scheme is said to be stable. To show that the scheme is stable it is then sufficient to show that, using the $l_1$-norm, 
	\begin{equation}
	\label{eq_stability}
		\sum_{i}\vert U(i,n)-u(i\Delta x,n\Delta t)\vert \leq C,
	\end{equation}
	for all $n$ where $C\in \mathbb{R}^+$. It is trivial to see that,
	\begin{equation}
		\sum_{i}\vert U(i,n)-u(i\Delta x,n\Delta t)\vert \leq \sum_{i}\vert U(i,n)\vert +  \sum_{i} \vert u(i\Delta x,n\Delta t)\vert.
	\end{equation}	
	From this, provided that the exact solution, $u(x,t)$, is bounded, the stability of the numerical method is assured if the approximation $U(i,n)$ remains bounded. 
	
	As $U(i,n)$ is derived from an underlying stochastic process and as such is a density we know that $U(i,n)\geq 0\;\forall i$ and hence,
	\begin{equation}
	\sum_{i}\vert U(i,n)\vert=\sum_{i} U(i,n)
	\end{equation}

	In general stability needs to be shown individually for the problem under consideration. However there are two general classes of problems where we can show unconditional stability. The first class we can consider is a periodic domain with a commensurate periodic potential. In such a case the boundary points are set to the interior point at the opposite boundary. Considering the case of the domain $x\in [0,l]$ and taking a lattice of $M+1$ points with boundary points at $i=0$ and $i=M$, then 
	\begin{equation}
	\label{eq_mass}
	\begin{split}
	\sum_{i=1}^{M}U(i,n)= &\sum_{i=1}^{M}U(i,n-1)+P_r(0,n-1)U(0,n-1)-P_l(1,n-1)U(1,n-1)\\
	&+P_l(M,n-1)U(M,n-1)-P_r(M-1,n-1)U(M-1,n-1).
	\end{split}
	\end{equation}
The periodic boundary conditions then imply that $U(0,n)=U(M-1,n)$, and $U(M,n)=U(1,n)$. The periodic potential ensures that the jump probabilities follow the same relation, $P_r(0,n)=P_r(M-1,n)$, and $P_l(M,n)=P_l(1,n)$. Thus the total mass is conserved as $\sum_{i=1}^{M}U(i,n)=\sum_{i=1}^{M}U(i,n-1)$ from Eq. (\ref{eq_mass}). As the initial mass is finite the conservation of mass guarantees the stability of the scheme.
	
	The next class of unconditionally stable problems is the case of zero-flux boundary conditions. Here we have the condition that no mass is entering the system at the boundary. For this to be true the continuum equations must satisfy the following conditions,
	\begin{equation}
	\label{eq_zfbc}
	\frac{\partial u(x,t)}{\partial x}\Big\vert_{x=0}= 2\beta F(0,t) u(0,t),
	\end{equation}
and	
	\begin{equation}
	\frac{\partial u(x,t)}{\partial x}\Big\vert_{x=l}= 2\beta F(l,t) u(l,t).
	\end{equation}
	Considering the case of the domain $x\in [0,l]$ and taking lattice of $M$ points the boundary conditions are enforced by setting ghost points at $i=0$ and $i=M+1$, then 
	\begin{equation}
	\label{eq_mass2}
	\begin{split}
	\sum_{i=1}^{M}U(i,n)= &\sum_{i=1}^{M}U(i,n-1)+P_r(0,n-1)U(0,n-1)-P_l(1,n-1)U(1,n-1)\\
	&+P_l(M+1,n-1)U(M+1,n-1)-P_r(M,n-1)U(M,n-1).
	\end{split}
	\end{equation}
Zero-flux boundary conditions are implemented by setting the ghost points such that the total mass will be conserved,
	\begin{align}
	\label{eq_nzfbc}
	U(0,n)&=\frac{P_l(1,n)}{P_r(0,n)}U(1,n),\\
	U(M+1,n)&=\frac{P_r(M,n)}{P_l(M+1,n)}U(M,n).
	\end{align}
Using the two point quadrature Boltzmann weights, Eq. (\ref{eq_boltz3}) and taking the limit $\Delta x \to 0$, Eq. (\ref{eq_nzfbc}) recovers Eq. (\ref{eq_zfbc}) as required. Again as the total mass is conserved the numerical scheme will be stable.

	\subsection{Accuracy}
	\label{sec_acc}
		In order to estimate the accuracy of the DTRW numerical scheme we consider the scheme as if it were derived from a Taylor series.
		Considering a non-linear advection diffusion equation of the form,
		\begin{equation}
		\frac{\partial u(x,t)}{\partial t}=D \frac{\partial^2 u(x,t)}{\partial x^2}-2\beta D \frac{\partial F(x,t,u) u(x,t)}{\partial x}.
		\end{equation}
		 The numerical scheme that approximates this equation is then,
		\begin{equation}\label{eq:DTRWscheme}
			u(x,t+\Delta t) =  p_l(x+\Delta x,t) u(x+\Delta x,t) + p_r(x-\Delta x,t) u(x-\Delta x,t),
		\end{equation}
		where, as before, the jump probabilities are given by Boltzmann weights, 
		\begin{equation}
			p_r(x,t) = \frac{\exp\left({\beta \frac{\Delta x}{2} (F(x+\Delta x,t)+F(x,t))}\right)}{\exp{\left(\beta \frac{\Delta x}{2} (F(x+\Delta x,t)+F(x,t))\right)}+\exp{\left(-\beta \frac{\Delta x}{2} (F(x-\Delta x,t)+F(x,t))\right)}},
		\end{equation}
		and 
		\begin{equation}
			p_l(x,t) = \frac{\exp\left({-\beta \frac{\Delta x}{2} (F(x-\Delta x,t)+F(x,t))}\right)}{\exp{\left(\beta \frac{\Delta x}{2} (F(x+\Delta x,t)+F(x,t))\right)}+\exp{\left(-\beta \frac{\Delta x}{2} (F(x-\Delta x,t)+F(x,t))\right)}}.
		\end{equation}
		Taking the Taylor expansion with respect to $\Delta x$, of the terms on the right hand side of Eq. (\ref{eq:DTRWscheme}), around $\Delta x=0$, to order $\mathcal{O}(\Delta x^3)$ we have
		\begin{equation}\label{eq:prexp}
		\begin{split}
			p_r(x-\Delta x,t)u(x-\Delta x,t) = &\frac{1}{2}u(x,t) +\frac{1}{2}\left( \beta F(x,t) u(x,t)- \frac{\partial u(x,t)}{\partial x}\right)\Delta x\\
			&- \frac{1}{4}\left(2 \beta u(x,t)\frac{\partial F(x,t)}{\partial x}+2 \beta F(x,t)\frac{\partial u(x,t)}{\partial x}-\frac{\partial^2 u(x,t)}{\partial x^2}\right) \Delta x^2 \\
			&-\frac{1}{24}\left(4\beta^3F(x,t)^3u(x,t)-12\beta \frac{\partial F(x,t)}{\partial x}\frac{\partial u(x,t)}{\partial x}\right.\\
			&\left.-9 \beta u(x,t)\frac{\partial^2 F(x,t)}{\partial x^2}-6\beta F(x,t)\frac{\partial^2 u(x,t)}{\partial x^2}+2\frac{\partial^3 u(x,t)}{\partial x^3}\right)\Delta x^3+ \mathcal{O}(\Delta x^4),
		\end{split}
		\end{equation}
		and
		\begin{equation}\label{eq:plexp}
		\begin{split}
			p_l(x+\Delta x,t)u(x+\Delta x,t) = &\frac{1}{2}u(x,t) -\frac{1}{2}\left( \beta F(x,t) u(x,t)- \frac{\partial u(x,t)}{\partial x}\right)\Delta x\\
			&- \frac{1}{4}\left(2 \beta u(x,t)\frac{\partial F(x,t)}{\partial x}+2 \beta F(x,t)\frac{\partial u(x,t)}{\partial x}-\frac{\partial^2 u(x,t)}{\partial x^2}\right) \Delta x^2 \\
			&+\frac{1}{24}\left(4\beta^3F(x,t)^3u(x,t)-12\beta \frac{\partial F(x,t)}{\partial x}\frac{\partial u(x,t)}{\partial x}\right.\\
			&\left.-9 \beta u(x,t)\frac{\partial^2 F(x,t)}{\partial x^2}-6\beta F(x,t)\frac{\partial^2 u(x,t)}{\partial x^2}+2\frac{\partial^3 u(x,t)}{\partial x^3}\right)\Delta x^3+ \mathcal{O}(\Delta x^4).
		\end{split}
		\end{equation}
		Substituting equations (\ref{eq:prexp}) and (\ref{eq:plexp}) into equation (\ref{eq:DTRWscheme}), dividing through by $\Delta t$ and collecting term appropriately we obtain,
		\begin{equation}\label{eq:DTRWacc}
			\frac{u(x,t+\Delta t) - u(x, t)}{\Delta t} = \frac{\Delta x^2}{2 \Delta t} \frac{\partial^2 u(x,t)}{\partial x^2} - \beta \frac{\Delta x^2}{\Delta t}  \frac{\partial F(x,t)u(x,t)}{\partial x}+ \mathcal{O}(\Delta x^4).
		\end{equation}
		Taking the Taylor expansion of the left hand side with respect to $\Delta t$ we have,
		\begin{equation}
			\frac{\partial u(x,t)}{\partial t}+\mathcal{O}(\Delta t) = \frac{\Delta x^2}{2 \Delta t} \frac{\partial^2 u(x,t)}{\partial x^2} - \beta \frac{\Delta x^2}{ \Delta t}  \frac{\partial F(x,t)u(x,t)}{\partial x}+ \mathcal{O}(\Delta x^4),
		\end{equation}
		with $D = \Delta x^2 / 2 \Delta t$ this reduces to,
		\begin{equation}
			\frac{\partial u(x,t)}{\partial t} = D \frac{\partial^2 u(x,t)}{\partial x^2} - 2 \beta D \frac{\partial F(x,t)u(x,t)}{\partial x}+\mathcal{O}(\Delta t) + \mathcal{O}( \Delta x^4).
		\end{equation}
		Finally noting that $\mathcal{O}(\Delta t) =\mathcal{O}(\Delta x^2) $ we see that the scheme will be convergent with order $\Delta x^2$.
		
\section{Burgers' Equation}
To show the utility of the numerical scheme we consider the viscous Burgers equation,
	\begin{equation}
	\label{eq_Burgers}
		\frac{\partial u(x,t)}{\partial t}=\nu \frac{\partial^2 u(x,t)}{\partial x^2}-u(x,t)\frac{\partial u(x,t)}{\partial x}.
	\end{equation}

The Burgers' equation may be transformed into the heat equation by using the Hopf-Cole transformation \cite{Hopf1950,Cole1951}, 
\begin{equation}
u(x,t)=-2 \nu \frac{1}{\phi(x,t)}\frac{\partial \phi(x,t)}{\partial x}.
\end{equation}
This transforms Burgers equation into the standard diffusion equation,
 \begin{equation}
 \frac{\partial \phi(x,t)}{\partial t}=\nu \frac{\partial^2 \phi(x,t)}{\partial x^2},
 \end{equation}
This allows for the construction of an infinite domain solution to the Burgers equation. One such solution is. 
		\begin{equation}
			u(x,t) = \frac{C_2+2\nu C_1^2 \tanh\left( t C_2-x C_1+C_3\right)}{C_1},
		\end{equation}
		where $C_1,C_2,C_3 \in \mathbb{R}$ are constants that depend on the initial condition. If we take the initial condition to be,
		\begin{equation}
		\label{eq_BInt}
		u(x,0)=1+2 \nu\tanh\left(c-x\right). 
		\end{equation}
		Then the solution is given by,
		\begin{equation}
		\label{eq_BSol}
		u(x,t)=1+2 \nu \tanh\left(c+t-x\right). 
		\end{equation}
		
	
To construct the numerical scheme we match the PDE of interest, Eq (\ref{eq_Burgers}), to the diffusive limit of the master equation, Eq. (\ref{eq_genAD}). 
	From this we see that,
	\begin{align}
	D&=\nu\\
	\beta&=\frac{1}{4\nu}\\
	F(x,t)&=u(x,t)
	\end{align}
	Thus the numerical scheme, Eq. (\ref{eq_NS}), becomes,
	 \begin{equation}
\label{eq_BNS}
\begin{split}
U(i,n)= &\frac{U(i-1,n-1)}{1+\exp(-\frac{\Delta x}{8\nu} \left(U(i-2,n-1)+2U(i-1,n-1)+U(i,n-1)\right))}\\
&+\frac{U(i+1,n-1)}{1+\exp(\frac{\Delta x}{8\nu} \left(U(i+2,n-1)+2U(i+1,n-1)+U(i,n-1)\right))}.
\end{split}
\end{equation}

\subsection{CFL Condition}
	
		The Courant-Friedrichs-Lewy (CFL) condition is a well known necessary criterion for the stability of finite difference type numerical schemes for hyperbolic partial differential equations. The CFL condition can be understood to imply that the numerical domain of the solution must be able to support the analytic solution of the PDE, otherwise information will be lost as the numerical solution is propagated in time \cite{Laney1998}. For the Burgers' equation under consideration CFL condition states that 
		\begin{equation}
			\frac{v \Delta t}{\Delta x} \le C_{\texttt{max}},
		\end{equation}
		where $v$ is the magnitude of velocity and $C_{\texttt{max}}$ is usually 1. We note here that the maximum velocity that our grid can accurately capture is $\Delta x / \Delta t$ and the maximum advective velocity prescribed by the model is $u(x,t)$. We then require 
		\begin{equation}
			\frac{\Delta x}{\Delta t} \le u(x,t) \quad \forall x,t,
		\end{equation}
		which is equivalent to the CFL condition in this case. Moreover, the relation
		\begin{equation}
			\Delta t = \frac{\Delta x^2}{2 D},
		\end{equation}
		indicates the relation between mesh spacing and the diffusive coefficient $D$, i.e. we have
		\begin{equation}
			\frac{\Delta x}{\Delta t} = \frac{2 D}{\Delta x} \le u(x,t).
		\end{equation}
		Unfortunately to impose this condition the solution, $u(x,t)$, needs to be known or estimated for all time and space. However, we may derive a useful heuristic from a physical argument: since the CFL condition ensures that the granularity of the mesh is sufficiently fine to accurately capture the velocity of propagation we may also interpret this as requiring the probability of jumping left or right to be within $(0,1)$, where a probability of jumping left (or right) of 1 (or 0) implies that a particle is moving at the maximum velocity allowed by the mesh, and within this framework this velocity may never be exceeded. We may then heuristically check that the jump probabilities remain sufficiently far from either 0 or 1. 
		
		We also argue that due to the viscosity present in the model, the maximum value of $u(x,t)$ occurs when $t=0$ for some value of $x$. We may then prescribe the CFL number to be $C_{\textrm{max}}=\smash{\displaystyle\max_{x}} \textrm{ } u(x,0)$.
		
		We note here that, in general, violation of the CFL condition does not guarantee that a numerical solution will blow up. We observe numerically that if the CFL condition is violated erroneous features form in the solution but the stability of the scheme, in respect to Eq. (\ref{eq_stability}) is not compromised.

	\subsection{Example 1: Dirichlet Boundary Conditions}
			
		 The infinite domain solution, Eq (\ref{eq_BSol}), may be used to find an exact solution of the Burger's equation on a finite domain, subject to time dependent Dirichlet boundary conditions. Taking Eq. (\ref{eq_Burgers}) on the domain $x\in [0,100]$ with the initial condition given by Eq. (\ref{eq_BInt}), and boundary conditions, 
		\begin{align}
			u(0,t)&=1+2 \nu \tanh (t+c),\\
			u(100,t)&=1+2 \nu \tanh (t+c-100).
		\end{align}
		The exact solution is then given by Eq. (\ref{eq_BSol}).
		
		The numerical scheme then comprises of Eq. (\ref{eq_BNS}), subject to the initial conditions,
		\begin{equation}
		U(i,0)=1+2 \nu\tanh\left(c-i \Delta x\right),
		\end{equation}
		The Dirichlet boundary conditions are implemented via Eq. (\ref{eq_DBC}), this requires the setting of the boundary points,
		\begin{align}
			U(0,n)&=1+2 \nu\tanh\left(n \Delta t+c\right)\\
			U(L,n)&=1+2 \nu\tanh\left(n \Delta t-100+c\right)
		\end{align}
		where is $L=\left\lfloor \frac{100}{\Delta x}\right\rfloor$. The jump probabilities for the boundary points must be computed using the single point quadrature form, Eq. (\ref{eq_boltz1pa}), so as to avoid the need for an additional point exterior to the boundary. Thus,
		\begin{equation}
			U(1,n)=\frac{U(0,n-1)}{1+\exp(-\frac{\Delta x}{2\nu}U(0,n-1))}+\frac{U(2,n-1)}{1+\exp(\frac{\Delta x}{8 \nu}(U(3,n-1)+2 U(2,n-1)+U(1,n-1)))},
		\end{equation}
	and
		\begin{equation}
			U(L-1,n)=\frac{U(L-2,n-1)}{1+\exp(-\frac{\Delta x}{8 \nu}(U(L-3,n-1)+2 U(L-2,n-1)+U(L-1,n-1)))}+\frac{U(L,n-1)}{1+\exp(\frac{\Delta x}{2 \nu}U(L,n-1))}
		\end{equation}
	For $2\leq i \leq L-2$, $U(i,n)$ will be given by Eq. (\ref{eq_BNS}).  
		
		A numerical solution was obtained where $c=-3$, and $\nu=0.45$. A range of $\Delta x$ values were chosen so that the numerical method would produce a value at $t=\frac{6250}{81}$, specifically $\Delta x \in \{\frac{25}{3},\frac{25}{12},\frac{25}{27},\frac{25}{48},\frac{1}{3},\frac{25}{108},\frac{25}{147},\frac{25}{192},\frac{25}{243},\frac{1}{12}\}$.  Confirming the accuracy analysis in section \ref{sec_acc} an accuracy of order $\Delta x^2$ was observed. Figure 1 provides a plot of the L$_1$-Norm of the error defined as,
		\begin{equation}
		E(n \Delta t)=\Delta x \sum_{i}\vert U(i,n)-u(i \Delta x, n\Delta t)\vert.
		\end{equation}
		
		\begin{figure}[h]
\begin{center}
\includegraphics[width=\textwidth]{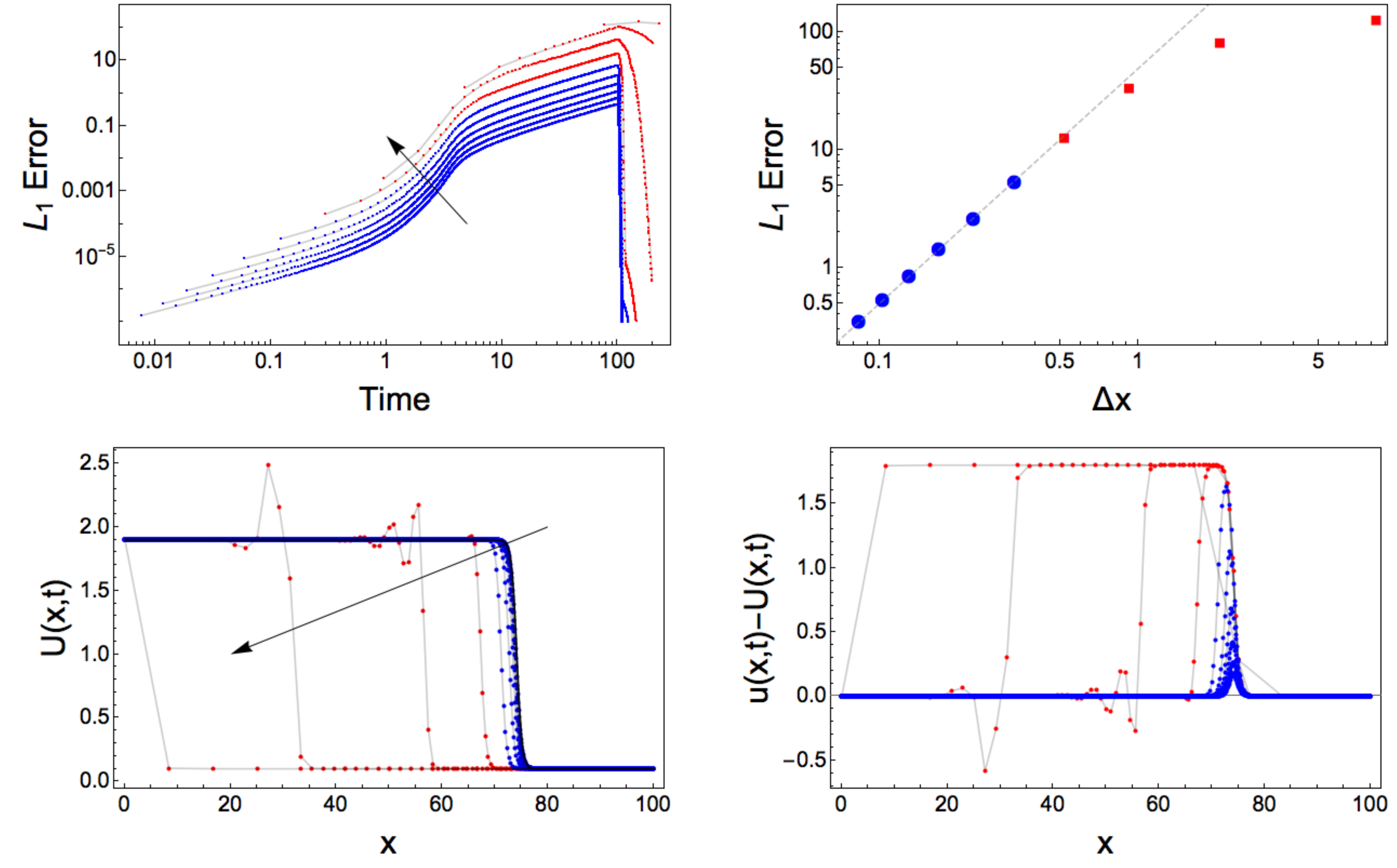}
\caption{Example One: Burgers Equation with time dependent Dirichlet boundary conditions. (Top Left) Error as a function of time. The arrow indicates direction of increasing $\Delta x$. Red points indicate that the CFL condition is broken.  (Top Right) Error as a function of $\Delta x$ at $t=\frac{6250}{81}$. (Bottom Left) Plot of the solutions at $t=\frac{6250}{81}$, the arrow indicates the direction of increasing $\Delta x$. The exact solution is shown as a solid black line. (Bottom Right) Difference between the exact solution and the numerical approximation at $t=\frac{6250}{81}$.}
\label{fig1}
\end{center}
\end{figure}

	\subsection{Example 2: Neumann Boundary Conditions}
		The second example we consider uses Neumann boundary conditions, again obtaining an exact solution from the infinite domain solution. Considering the Burger's equation, Eq. (\ref{eq_Burgers}), on the domain, $x\in[0,100]$ with the initial condition given by Eq. (\ref{eq_BInt}), and with time dependent Neumann boundary conditions given by,
		\begin{align}
		\frac{\partial u(x,t)}{\partial x}\big\vert_{x=0}&=-2\nu\; \mathrm{sech}^2(t+c),\\
		\frac{\partial u(x,t)}{\partial x}\big\vert_{x=10}&=-2\nu\; \mathrm{sech}^2(t-100+c).
		\end{align}
Then the solution will once again be given by Eq. (\ref{eq_BSol}). In this case we need the boundary points, $x=0$ and $x=100$, to lie between lattice points. As such we have,
	\begin{equation}
		U(i,n)\approx u((i-\frac{1}{2})\Delta x,n \Delta t).
	\end{equation}
Similarly to Example 1, the numerical scheme comprises of Eq. (\ref{eq_BNS}), subject to the initial conditions,
		\begin{equation}
		U(i,0)=1+2 \nu\tanh\left(c-(i-\frac{1}{2}) \Delta x\right),
		\end{equation}
The boundary conditions are now implemented by setting ghost points outside the domain via Eq. (\ref{eq_NBC}), so that,
	\begin{equation}
		U(L+1,n)=U(L,n) -2\Delta x \nu\; \mathrm{sech}^2(n \Delta t-100+c),
	\end{equation}
and
	\begin{equation}
		U(0,n)=U(1,n)+2 \Delta x \nu\; \mathrm{sech}^2(n \Delta t+c).
	\end{equation}

Once again, the jump probabilities for the ghost points must be computed using the single point quadrature form, Eqs. (\ref{eq_boltz1p}) and (\ref{eq_boltz1pa}), so as to avoid the need for an additional ghost point.  

Similarly to example one, values of $\Delta x$ were chosen so that the solution could be compared at $t=\frac{6250}{81}$, specifically $\Delta x \in \{\frac{25}{3},\frac{25}{12},\frac{25}{27},\frac{25}{48},\frac{1}{3},\frac{25}{108},\frac{25}{147},\frac{25}{192},\frac{25}{243},\frac{1}{12}\}$. As in example one, an accuracy of order $\Delta x^2$ was observed. 

	\begin{figure}[h]
		\label{figEx2}
		\begin{center}
			\includegraphics[width=\textwidth]{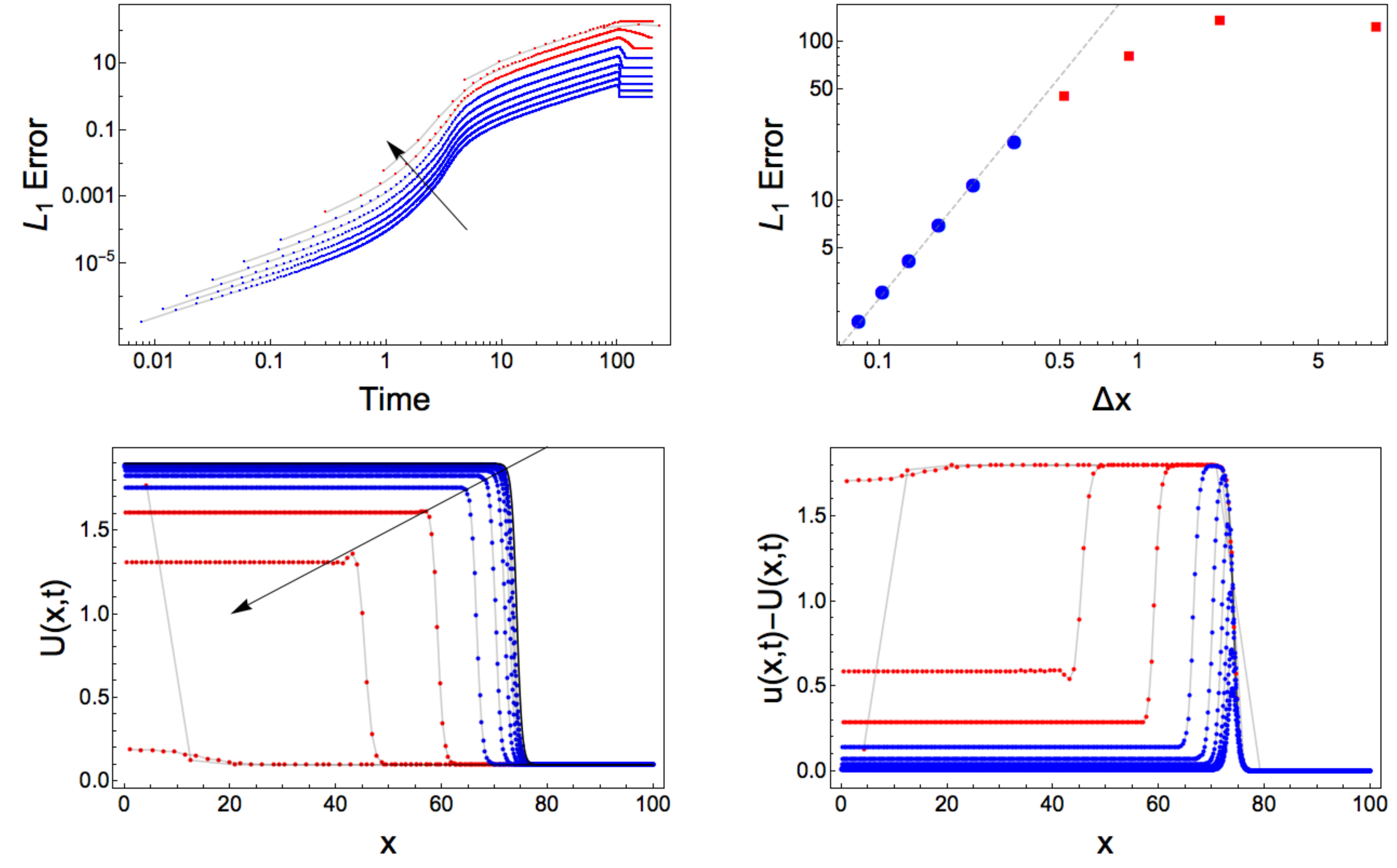}
			\caption{Example Two: Burgers Equation with time dependent Neumann boundary conditions. (Top Left) Error as a function of time. The arrow indicates direction of increasing $\Delta x$. Red points indicate that the CFL condition is broken.  (Top Right) Error as a function of $\Delta x$ at $t=\frac{6250}{81}$. (Bottom Left) Plot of the solutions at $t=\frac{6250}{81}$, the arrow indicates the direction of increasing $\Delta x$. The exact solution is shown as a solid black line. (Bottom Right) Difference between the exact solution and the numerical approximation at $t=\frac{6250}{81}$.}
		\label{fig2}
		\end{center}
	\end{figure}
		
\section{Conclusion}

	We have shown that a discrete time random walk (DTRW) can be used to construct a numerical method for the solution of general non-linear advection-diffusion partial differential equations. The derived explicit numerical scheme is both easy to implement and stable, even in the case where the underlying equations admit shock solutions. Provided that the underlying stochastic process that the numerical scheme is based on exists, we have shown that the scheme must be stable.  As a prototypical example of a nonlinear advection-diffusion equation we have considered the one-dimensional viscous Burgers' equation. Burgers' equation arises from the simplification of the homogenous incompressible Navier-Stokes equations and has been extensively used as a benchmark for numerical methods for hyperbolic partial differential equations due to the non-linearities and array of dynamics exhibited by the equation. The presented DTRW scheme accurately captures the dynamics of Burgers' equation when compared with exact solutions.

\section*{Acknowledgments}
This work was supported by the Australian Commonwealth Government (ARC DP140101193). B.A.J. acknowledges the National Research Foundation of South Africa under grant number 94005. Opinions expressed and conclusions arrived at are those of the author and are not necessarily to be attributed to the CoE-MaSS.

\bibliographystyle{elsarticle-num}
\bibliography{DTRW}
	
\end{document}